\documentclass[opre,nonblindrev]{informs3} 

\DoubleSpacedXII 

\usepackage[numbers]{natbib}
 \bibpunct{[}{]}{,}{n}{}{,}%
 \def\bibfont{\small}%
 \makeatletter
\renewcommand\@biblabel[1]{#1.}
\makeatother

\def\b{\color{black}}

\TheoremsNumberedThrough     
\ECRepeatTheorems

\EquationsNumberedThrough    

\usepackage{graphicx, subfigure}
\usepackage{amsfonts,mathtools}
\usepackage{verbatim}
\usepackage{txfonts}
\usepackage{cases}
\usepackage{color}
\usepackage{url}
\usepackage[normalem]{ulem}

\begin{document}


\RUNAUTHOR{E. Feinberg and X. Zhang}

\RUNTITLE{
Optimal Control of $M/M/\infty$
}

\TITLE{\uppercase{Optimal Switching On and Off the Entire Service Capacity of a Parallel Queue}}

\ARTICLEAUTHORS{%
\AUTHOR{{\large E}ugene {\large F}einberg}
\AFF{Department of Applied Mathematics and
Statistics, Stony Brook University, \\
Stony Brook, NY 11794,
\EMAIL{efeinberg@notes.cc.sunysb.edu}}
\AUTHOR{{\large X}iaoxuan {\large Z}hang}
\AFF{IBM T.J. Watson Research Center,\\
 Yorktown Heights, NY 10598,
 \EMAIL{zhangxiaoxuan@live.com}}
} 

\ABSTRACT{ This paper studies optimal switching on and off of the
entire service capacity of an $M/M/\infty$ queue with holding,
running and switching costs where the running costs depend only on
whether the system is  running or not. The goal is to minimize average
costs per unit time. The main result is that an average-optimal
policy either always runs the system or is an $(M,N)$-policy
defined by two thresholds $M$ and $N$, such that the system is
switched on upon an arrival epoch when the system size accumulates
to $N$ and  is switched off upon a departure epoch when the system
size decreases to $M$.  It is shown  that this optimization
problem can be reduced to a problem with a finite number of
states and actions, and an average-optimal policy can be computed
via linear programming. An example, in which the optimal
$(M,N)$-policy outperforms the best $(0,N)$-policy, is provided.
Thus, unlike the case of single-server queues studied in the
literature, $(0,N)$-policies may not be average-optimal.}



\maketitle

%


\section{\uppercase{Introduction}}\label{secMF}

This paper studies 
optimal control of a parallel M/M/$\infty$
queue with Poisson arrivals  and an unlimited number of
independent identical servers with exponentially distributed
service times.
The
cost to switch the system on is {\b $s_1$} and the cost to switch the
system off is {\b $s_{0}$}.  The other costs include the linear holding
cost $h$ for {\b each} unit of time that {\b a} customer spends in {\b the}
system, the running cost $c_1$ per unit time when the system is on
and the idling cost $c_0$ per unit time when the system is off. It
is assumed that $s_0,s_1\ge 0$, $s_0+s_1>0$, $h>0$, and $c_1>c_0.$
Denote $c=c_1-c_0$. Without loss of generality, let $c_0=0$ and
$c_1=c>0.$  The goal is to {\b minimize}  average costs per unit time.

The main result of this paper is that either the policy that
always keeps the system on is average-optimal or, for some integers
$M$ and $N$, where $N>M\ge 0$, the so-called $(M,N)$-policy is
average-optimal. The $(M,N)$-policy switches the running system
off when the number of customers in the system is not greater than
$M$ and it switches the idling system on when the number of customers
in the queue reaches or exceeds $N.$  { It is shown in this paper
that this optimization problem can be reduced to a problem with
finite number of states and actions and an average-optimal policy
can be computed via linear programming.  An example when the best
$(0,N)$-policy is not optimal is provided.}
%


Studies on control problems for queues started around fifty years
ago, and one of the first papers on this topic,
\citet{yadin1963queueing}, dealt with switching on and off a
server of a single-server queue. 
 \citet{heyman1968} 
showed the optimality of {\b a} $(0,N)$-policy, which is usually
called
an $N$-policy, for $M/G/1$ queues.  
%
  The early results on switching servers in
single-server queues led to two relevant research directions:

(i) optimality of $(0,N)$-policies or their ramifications under
very general assumptions such as batch arrivals, start-up and
shut-down times and costs, nonlinear holding costs, known
workload, and so on; see \citet{lee1989control,
federgruen1991optimality, altman1993optimal,
denardo1997stochastic}{\b, and} \citet{ feinberg2002optimality};

 (ii) decomposition results for queues with vacations; see
\citet{fuhrmann1985stochastic, hofri1986queueing,
shanthikumar1988stochastic,kella1989threshold}, and
\citet{kella1991queues}. \citet{sobel1969} studied
$(M,N)$-policies for $GI/G/1$ queues.


{\b As for general multi-server parallel queues, switching on and
off individual servers for a parallel queue is a  more
difficult problem. 
 Even for an $M/M/n$ queue,  there's no known description of an optimal switching policy
for individual servers when $n>2$; see
\citet{bell1975technical,bell1980optimal},
\citet{rhee1990distribution}, and \citet{igaki1992exponential}.
%
%
Studies of stationary distributions and performance evaluation{\b
s} {\b for} parallel queues with vacations  (\citet{levy1976m},
\citet{huang1977optimal}, \citet{kao1991analyses},
\citet{browne1995parallel}, \citet{chao1998analysis} and
\citet{li2000optimal}) usually assume that vacations begin when
the system is empty. Observe that, if vacations start when the
system becomes empty and end simultaneously for all the servers,
the model describes a particular case of switching the entire
service capacity of the system on and off.
\citet{browne1995parallel} studied $M/G/\infty$ {\b queues with
vacations} and described how to compute  the best $(0,N)$-policy
for switching on and off the entire service capacity. 


{\b This research is motivated by two  observations: (i) the
problem of switching  on and off the entire service capacity of 
the facility
has an explicit solution described in this paper while there is no 
known explicit solution for problems with 
  servers that can be switched on and off individually, and
 (ii) with the development of internet and high performance
 computing, many applications behave in the way described in this paper.  For example, consider a  a service provider that uses cloud computing and pays for the
 time the cloud is used.  When there are many service
 requests, it {\b is worth paying} for using the cloud, and when there is a small
 number of service requests, it may be too expensive to use the
 cloud.  This paper analyzes such a situation and finds an optimal solution.  In fact,
 the performance analysis literature treats 
 cloud computing as a
parallel queue with thousands of servers (see \citet{ibmscp}), and 
the number of servers in the models will increase with the
development of technologies.
{\b Many papers and research works on cloud computing 
model cloud computing  facilities as  multi-server queues; see 
\citet{mazzucco2010maximizing} and \citet{khazaei2012performance}.
\citet{mazzucco2010maximizing} {\b{studies} } the revenue
management problem from 
the perspective of a cloud computing provider and investigates 
the resource allocation via dynamically powering the servers on or
off. }
However, there can be a huge number of servers in a cloud
computing center, typically of the order of hundreds or thousands;
see \citet{greenberg2008cost, msft}, and \citet{amazon}. Given
that the number of servers is large and tends to increase over
time with the development of new technologies, it is natural to
model controlling of the facility as an $M/M/\infty$ queue rather
than an $M/M/n$ queue if this leads to analytical advantages.
Here we study a model based on an M/M/$\infty$ queue and find a
solution.

In addition to cloud computing, another {\b example} comes from the
application to IT software maintenance.
\citet{kulkarni2009optimal} studied the software maintenance
problem as a control problem for a queue  formed by software
maintenance requests { generated by} software bugs  experienced by
customers. Once a customer is served and the appropriate bug is
fixed in the new software release or patch, it also provides
{\b solutions} to some other customers in the queue {\b and these customers are served simultaneously}.  In
\citet{kulkarni2009optimal}, it was assumed that the number of
customers leaving the queue at a service completion time has a
binomial distribution. This problem was modeled in
\citet{kulkarni2009optimal} as an optimal switching problem for an
$M/G/1$ queue in which a binomially distributed number of
customers depending on the queue size {\b are} served each time, and the
best policy was found among the policies that switch the system off
when it is empty and switch it on when there are $N$ or more
customers in the system. Here we observe that after an appropriate
scaling, the software maintenance  problem with exponential
service times and the optimal switching problem for an
$M/M/\infty$ queue have the same fluid approximations.  So, the
result on average-optimality of $(M, N)$-policies described {\b here}
provides certain insights to the {\b software maintenance} problem studied in
\citet{kulkarni2009optimal}.

{\b As a conclusion} to the introduction, we describe the
structure of the paper and some of the main ideas.  There are two
main obstacles in the analysis of the M/M/$\infty$ switching
problem compared to a 
 single-server one.  First, the service
intensities are  unbounded, and therefore the standard reduction
of continuous-time problems to discrete time via uniformization
{\b can not be applied}. Second, there are
significantly more 
known decomposition and performance analysis results for
single-server queues than for parallel queues and, in particular,
we are not aware of such results for M/M/$\infty$ queues with
vacations that can start when the queue is not empty.  The first
obstacle is resolved  by reducing the discounted version of the
problem to negative dynamic programming instead of to 
discounted dynamic programming. The second obstacle {\b is}
resolved by solving a discounted problem {\b for the system that
cannot be switched off}. This problem is solved by using optimal
stopping, where the stopping decision corresponds to starting the
servers, and its solution is used to derive useful inequalities
and to reduce the problem for the original} M/M/$\infty$ queue to
a control {\b problem} of a semi-Markov process with finite state
and action sets representing 
 the system being always on when the
number of customers exceeds a certain level.

The optimal switching problem for an M/M/$\infty$ queue is modeled
in Section~\ref{secPF} as a Continuous-Time Markov Decision
Process (CTMDP) with unbounded transition rates.  Such a CTMDP
cannot be reduced to discrete time via uniformization;
{\b see}, e.g., \citet{piunovskiy2012transformation}.  The available results for
average costs require that any stationary policy defines an
ergodic continuous-time Markov chain; see \citet{guo2002dsc} and
\citet[Assumptions 7.4. 7.5 on p. 107]{guodhern2009}. These
assumptions do not hold for the problem we consider because the
policy that always keeps the system off defines a transient Markov
chain. Therefore, in this paper we provide a direct analysis of
{\b the} problem.

{\b Section \ref{secDisc} } studies  expected total discounted
costs.  Such a CTMDP can be reduced to a discrete-time problem
with the expected total costs; see \citet{feinberg2012reduction},
\citet{piunovskiy2012transformation}. Since transition rates are
unbounded,  expected total costs for the discrete-time problem 
cannot be presented as expected total discounted costs with the
discount factor smaller than 1.  However, since all the costs are
nonnegative, the resulting discrete-time problem belongs to the
class of negative MDPs {\b that deal with minimizing expected
total nonnegative costs, which is equivalent to maximizing
expected total nonpositive rewards}. For this negative MDP we
derive the optimality equation, show that the value function is
finite, and establish the existence of stationary discount-optimal
policies; see Theorem~\ref{Th0}.

Subsection~\ref{sec_full}  investigates the discounted total-cost 
problem limited to the policies that never {\b switch} the running
system off.  Such policies are called full-service {\b policies}.
By using the fact that the number of customers in an M/G/$\infty$
queue at each time has a Poisson distribution (see \citet[p.
70]{ross1996stochastic}), we compute {\b value functions} for
full-service policies and then compute the discount-optimal
full-service policy in Theorem~\ref{Th2}. This is done by
analyzing the optimality equation for an optimal stopping problem
when stopping, in-fact, corresponds to the decision to start the
system. The optimal full-service policy is defined by a number
$n_\alpha$ such that the system should be switched on as soon as
the number of customers is greater than or equal to $n_\alpha$,
where $\alpha>0$ is the discount rate.  The important feature of
the function $n_\alpha$ is that it is increasing in $\alpha$ and
therefore bounded when $\alpha\in(0,\alpha^*]$ for any $\alpha^*
\in (0,\infty).$ In Section~\ref{sec_red}, the problem with the
expected discounted total {\b costs} is reduced to a problem with
finite state and action sets by showing in
Lemma~\ref{lem_Mbdd_old} that the system should always be on, if
the number of customers is greater or equal than $n_\alpha.$
In Section~\ref{s:avg}, by using the vanishing discount rate
arguments, we prove the existence of stationary average-optimal
policies and describe  {\b the $(M,N)$-policy in Theorem \ref{t:avopt}}.
A linear program (LP) for their computation is provided in
Section~\ref{s:eg}. An example {\b in which} the best $(0,N)$-policy is not
average-optimal is given in Section~\ref{s:no}.

\section{\uppercase{Problem Formulation}}\label{secPF}
We model the above described control problem for an $M/M/\infty$
queue as a CTMDP with a
countable state space and a finite number of actions; see \citet{kitaev} and \citet{guodhern2009}. In general, such a CTMDP is defined by {\b the} tuple
$\{Z,A,A(z),q,c\},$ where $Z$ is a countable state space, $A$ is
a finite action set, $A(z)$ are sets of actions available in
state{\b s} $z\in Z$, and $q$ and $c$ are transition and cost rates,
respectively.  A general policy can be time-dependent,
history-dependent, and at a jump epoch the action that controls
the process  is the action selected at the previous state; see \citet[p.
138]{kitaev}. An initial state $z\in Z$ and a policy $\pi$ define
a stochastic process $z_t$ and {\b the } expectations for this stochastic
process are denoted by $E_z^\pi$. Let $C({\b t})$ be {\b the }cumulative costs
incurred during the time interval $[0,t].$ For $\alpha>0$, the
expected
total discounted cost is 
\begin{align}\label{eqObjEF}
&V_{\alpha}^{\pi}(z)=E^{\pi}_{z}\int_0^\infty e^{-\alpha t}dC(t),
\end{align}
and the average cost per unit time is 
\begin{align} \label{eqObjEF1}
v^{\pi}(z)=\limsup_{t\rightarrow\infty}t^{-1}E^{\pi}_{z}C(t).
\end{align}
Let
\begin{eqnarray}
&V_{\alpha}(z)=\displaystyle\inf_{\pi}V_{\alpha}^{\pi}(z),\label{optV}\\
&v(z)=\displaystyle\inf_{\pi}v^{\pi}(z).\label{optv}
\end{eqnarray}
A policy $\pi$ is called discount-optimal if
$V^{\pi}_{\alpha}(z)=V_{\alpha}(z)$  for all initial states $z\in
Z$. A policy $\pi$ is called average-optimal if
$v^{\pi}(z)= v(z)$ for all initial states $z\in Z$.

For our problem, the states of the system change only at  arrival
and departure epochs, which we call jump
epochs.  
The state of the system at time $t\ge 0$ is $z_t=(x_t,\delta_t),$
where $x_t$ is the number of customers in the system at time $t,$
and $\delta_t$ is the status of the servers that an arrival or
departure saw at the last jump epoch.  If $\delta_t=0$, the severs
at the last jump epoch during the interval $[0,t]$ were off, and,
if $\delta_t=1$, they were on. In particular, if the last jump
epoch was departure, $\delta_t=1.$  If the last jump epoch was an
arrival, then $\delta_t=1$, if the servers were on at that epoch,
and $\delta_t=0$ otherwise.  The initial state
$z_0=(x_0,\delta_0)$ is given.

The state space is $Z=\mathbb{N}\times\{0,1\}$, where
$\mathbb{N}=\{0,1,\ldots\}$, and the action set {\b is}
$A=\{0,1\}$, with 0 means that the system should be off and 1
meaning that the system should be on.  If at time $t$ the state is
$z_t=(x_t,\delta_t)$, this means that $x_t$ is the number of
customers in the system at time $t$, and $\delta_t\in \{0,1\}$ is
the control used at the last jump epoch during the interval
$[0,t]$.  The action sets $A(z)=A$ for all $z\in Z.$
A stationary policy chooses actions deterministically at jump epochs
and follows them until the next jump. In addition, the choice of
an action depends only on the state of the system $z=(x,\delta)$,
{\b where $x$ is the number of customers in the system and $\delta\in\{0,1\}$  is the status
of the system prior to the last jump}.

%
The transition rate from a state $z=(i,\delta)$ with {\b  an
}action $a\in A$ to a state $z'=(j,a)$, where $j\ne i,$ is
$q(z'|z,a)=q(j|i,a)$, with
\begin{align}\label{q}
q(j|i,a)
=
\begin{cases}
\lambda, &\text{if $ j=i+1$,}\\
i\mu{\b ,} &\text{if $i>0,\ a=1,\  j=i-1$,}\\
0, &\text{otherwise},
\end{cases}
\end{align}
{\b where $\lambda$ is the intensity of the arrival process and $\mu$ is the service rate of individual servers.}
At state $z=(i,\delta)$, define
$q(z,a)= q(i,a)=
\displaystyle{\sum_{j\in\mathbb{N}\setminus\{i\}}}q(j|i,a)$ and
$q(z|z,a)=q(i|i,a)=-q(i,a).$

%
The costs include the linear holding cost $h$  for a
unit time that a customer spends in the system, the running cost
$c$ per unit time {\b when} the system is on, the start-up cost $s_1$,
{\b and} the shut-down cost $s_0$, where $h,c>0$, $s_0,s_1\ge0$, and
$s_0+s_1>0.$  At state $z=(i,\delta)$, if action 1 is taken,   the
cost rate is $hi+c$; if action 0 is taken,  the cost rate is $hi$.
 At state $z=(i,\delta)$, if action 1 is
taken, the {\b instantaneous} cost  $(1-\delta)s_1$ is incurred;
if action 0 is taken, the {\b instantaneous} cost $\delta s_0$ is
incurred. The presence of {\b instantaneous} switching costs $s_0$
and $s_1$ complicates the situation, because standard models of
CTMDPs deal only with cost rates.
  To resolve
this complication, observe that, since $s_0+s_1>0$,  the number of
times  the system{\b 's} status (on or off) changes up to
any time {\b $N(t)$},  when $t<\infty$, should be finite with probability 1 for a
policy $\pi$ and an initial state $z$.  Otherwise,
$V_{\alpha}^{\pi}(z)=v^{\pi}(z)=\infty$ for all $\alpha>0.$ Let
$0\le t_1<t_2<\ldots$ be the times when the system is switched on
or off.  Let $a(t)$ be {\b the} action 0 or 1, selected at time t.
If this function has a finite number of jumps on each finite
interval, we consider the function $a$ being left continuous. This
is consistent with the definition of a general policy as a
predictable function; see \citet[p. 138]{kitaev}. {\b For} an
initial state $z$ and a policy $\pi$, {\b if} the value of $N(t)$
is finite with probability 1 for all finite $t$, we define
%
\begin{align*}
C(t)=\int_0^t\left(h
x_t+ca(t)\right)dt+\sum_{n=1}^{N(t)}s_{a(t_n+)}|a(t_n+)-a(t_n)|,
\end{align*}
and consider the expected total discounted costs $V_\alpha^\pi(z)$
and the expected average costs  per unit time $v^\pi(z)$ defined
in (\ref{eqObjEF}) and (\ref{eqObjEF1}), respectively. {\b For some $t<\infty$, } if
$N(t)=\infty$ {\b with positive probability}, we set
$V_\alpha^\pi(z)=v^\pi(z)=\infty.$  

\section{\uppercase{Discounted Cost Criterion}}\label{secDisc}

In this section we study the expected total {\b discounted} cost
criterion. {\b We first reduce the CTMDP to the discrete time MDP
and {\b then} study the so called full-service policies, which are
 used to reduce the original problem to an
equivalent problem with a finite state space. }

\subsection{\uppercase{Reduction to Discrete Time and Existence of Stationary Discount-Optimal
Policies}}\label{sec_OE}

In this subsection, we formulate the optimality equation and prove
the existence of stationary discount-{\b optimal policy}.  This
is done by {\b reducing the} problem to discrete time.

When the system is on and there are i customers, the time until
the next jump has an exponential distribution with intensity
$q(i,1)=\lambda+i\mu\to\infty$ as $i\to\infty.$ Since the jump
rates are {\b unbounded}, it is impossible to present the problem as a
discounted MDP {\b in discrete-time} with {\b a} discount {\b factor} smaller than 1.  Thus, we
shall present our problem as minimization of the expected total
costs. If the decisions are chosen only at jump times, the expected
{\b total discounted {\b sojourn} time} until the next jump epoch is
$
\tau_\alpha(z,a)
{\b=\int_0^\infty(\int_0^t e^{-\alpha s}ds)q(z,a)e^{-q(z,a)t}dt}
=\int_0^\infty e^{-\alpha
t}e^{-q(z,a)t}dt=\displaystyle\frac{1}{\alpha+q(z,a)}, $ and
the one-step cost is $
C_\alpha(z,a)=|a-\delta|s_a+(hi+ac)\tau_\alpha(z,a) $ with
$z=(i,\delta).$ For $\alpha=0$, we  denote   $\tau_{0}(z,a)$ and
$C_{0}(z,a)$ as $\tau(z,a)$ and $C(z,a)$  respectively.

By \citet[Theorem 5.6]{feinberg2012reduction}, there exists a
stationary discount-optimal policy, the value function
$V_\alpha(z)$ satisfies the discount-optimality equation,
\begin{equation} \label{eqF2}
V_\alpha(z)=\min_{a\in A(z)}\{C_{\alpha}(z,a)+\sum_{z^{\prime}\in
Z}\frac{q(z'|z,a)}{\alpha+q(z,a)}V_\alpha(z')\},\qquad
z\in Z{\b ,}
\end{equation}
and a stationary  policy $\phi$ is discount-optimal if and only if
\begin{equation} \label{eqF1} V_\alpha(z)=C_{\alpha}(z,\phi(z))+\sum_{z^{\prime}\in
Z}\frac{q(z'|z,{\b \phi(z)})}{\alpha+q(z,\phi(z))}V_\alpha(z),
\quad z\in Z.
\end{equation}
Formulae (\ref{eqF2}) and  (\ref{eqF1})  {\b imply} that the discounted
version of the problem is equivalent to finding a policy
{\b that minimizes} the expected total costs
 for {\b the} discrete-time MDP $\{Z,A,A(z),p_\alpha,C_\alpha\}$ with sub-stochastic
transition probabilities
$p_\alpha(z'|z,a)=q(z'|z,a)/\left(\alpha+q(z, {\b a})\right)$ and with
one-step cost $C_\alpha(z,a),$  where $\alpha>0.$

As mentioned above, classic CTMDPs do not deal with instantaneous
costs described in the previous section.  However, if we replace
the instantaneous costs $s_a$, {\b $a\in\{0,1\}$}, with the cost
rates $s(z,a)=s_a|a-\delta|(\alpha+q(z, a))$, where
$z=(i,\delta)$, then the expected total discounted cost until the
next jump does not change for policies that change actions only at
jump epochs. For an arbitrary policy, the expected total
discounted  cost until the next jump can either decrease or remain
unchanged, if instantaneous costs are replaced with the described
cost rates. This follows from Feinberg~\cite[Theorem
1]{feinberg1994optimality}, which implies that the defined cost
rates $s(z,a)$ correspond to the situation when only the first
nonzero switching cost after the last jump is charged and the
remaining switchings are free (in particular, if $s_0,s_1>0$, only
the first switching is charged). Thus, a discount-optimal policy
for the problem with the switching cost rates $s(z,a)$ is also
discount-optimal for the original {\b problem} with instantaneous
switching costs, and the optimality equation  \eqref{eqF2} is also
the optimality equation for the original problem with the goal to
minimize the expected total discounted costs.

{\b The following lemma computes the value function under the
policy that always runs the system.} {\b This lemma produces an
upper bound for the value function $V_\alpha$ and, in addition, it
shows that the value function takes finite values.} 

\begin{lemma}\label{lem_U}
Let $\phi$ be a policy that always runs the system.  For all
$i=0,1,\ldots$,
\begin{align}\label{U}
V^{\phi}_{\alpha}(i,\delta)=
(1-\delta)s_1+\frac{hi}{\mu+\alpha}+\frac{h\lambda}{\alpha(\mu+\alpha)}+\frac{c}{\alpha}<\infty.
\end{align}
\end{lemma}
\proof{Proof.} $V^\phi_\alpha(i,0)=s_1+V^\phi_\alpha(i,1)$,
 or equivalently,
$V^\phi_\alpha(i,\delta)=(1-\delta)s_1+V^\phi_\alpha(i,1).$
Observe that
\begin{align}\label{V1_full}
V^\phi_{\alpha}(0,1)=E\left[\int_0^{\infty}e^{-\alpha
t}\left(hX_0(t)+c\right)dt\right]
=hE\left[\int_0^{\infty}e^{-\alpha
t}X_0(t)dt\right]+\frac{c}{\alpha}=\frac{h\lambda}{\alpha(\mu+\alpha)}+\frac{c}{\alpha},
\end{align}
where $X_0(t)$ is the number of busy servers at time $t$ if at
time 0 the system is empty.  The last equality in (\ref{V1_full})
holds because, according to Page 70 in \citet{ross1996stochastic},
$X_0(t)$ has a Poisson distribution with the mean
$\lambda\int_0^te^{-\mu
t}dt=\displaystyle\frac{\lambda}{\mu}\left(1-e^{-\mu t}\right)$.
Thus,
\begin{align*}
E\left[\int_0^{\infty}e^{-\alpha
t}X_0(t)dt\right]=\int_0^{\infty}e^{-\alpha
t}\frac{\lambda}{\mu}\left(1-e^{-\mu t}\right)dt
=\frac{\lambda}{\alpha(\mu+\alpha)}.
\end{align*}
Also observe that
$$V^\phi_\alpha(i,1)=G_\alpha(i)+V^\phi_\alpha(0,1)=iG_\alpha(1)+V^\phi_\alpha(0,1),$$
where $G_\alpha(i)$ is the expected total discounted holding cost
to serve $i$ customers that are in the system at time 0.  Since
the service times are exponential, $G_\alpha(1)=E\left[\int_0^\xi
e^{-\alpha t}hdt\right]=\displaystyle\frac{h}{\mu+\alpha},$ where
$\xi\sim\exp(\mu).$ 
 In addition, $V_\alpha^\phi(i,0)=s_1+V_\alpha^\phi(i,1).$ \hfill\Halmos
\endproof

 We follow the conventions that
$V_\alpha(-1,\delta)=0,$
 $\displaystyle\sum_{\emptyset}=0,$ and $\displaystyle\prod_{\emptyset}=1.$
The following theorem is the main result of this subsection.

\begin{theorem}\label{Th0}  For any $\alpha>0$ the following statements
hold:    
\begin{hitemize}
\item[(i)] for all $i=0,1,\ldots,$
\begin{equation}\label{V1_phi0}
V_\alpha(i,\delta)\le
(1-\delta)s_1+\frac{hi}{\mu+\alpha}+\frac{h\lambda}{\alpha(\mu+\alpha)}+\frac{c}{\alpha};
\end{equation}
\item[(ii)]   for all $i=0,1,\ldots$ and {\b for }all
$\delta=0,1,$ the value function $V_{\alpha}(i,\delta)$ satisfies
the discount-optimality equation
\begin{align}\label{DCOE}
V_{\alpha}(i,\delta)
=  \min_{a\in\{0,1\}} &\left\{C_{\alpha}((i,\delta),a)+\frac{q(i-1|i,a)}{\alpha+q(i,a)}V_\alpha(i-1,a)+\frac{q(i+1|i,a)}{\alpha+q(i,a)}V_\alpha(i+1,a)\right\}\notag \\
=  \min &\left\{
(1-\delta)s_1+\displaystyle\frac{hi+c}{{\alpha+\lambda+i\mu }}
+\frac{\lambda}{\alpha+\lambda+i\mu }V_{\alpha}(i+1,1)
+\frac{i\mu}{\alpha+\lambda+i\mu }V_{\alpha}(i-1,1),\right.\notag\\
&\left.\ \ \ \delta
s_0+\frac{hi}{\alpha+\lambda}+\frac{\lambda}{\alpha+\lambda}V_{\alpha}(i+1,0)\right\};
\end{align}
 \item[(iii)] there exists a stationary discount-optimal policy, and
 a stationary policy $\phi$ is discount-optimal if and only if for all
 $i=0,1,\ldots$ and for all $\delta=0,1$,
 \begin{align*}
V_{\alpha}(i,\delta) =  \min_{\phi(i,\delta)\in\{0,1\}}
\left\{C_{\alpha}((i,\delta),\phi(i,\delta))+\frac{q(i-1|i,a)}{\alpha+q(i,\phi(i,\delta))}
V_\alpha(i-1,\phi(i,\delta)) \ \right.\\
\left.
+\frac{q(i+1|i,a)}{\alpha+q(i,\phi(i,\delta))}V_\alpha(i+1,\phi(i,\delta))\right\}.
\end{align*}
\end{hitemize}
\end{theorem}
\proof{Proof.}
(i) 
Consider the policy $\phi$ that always runs the system. 
 Then $V_\alpha(i,\delta)\le
V_\alpha^{\phi}(i,\delta),$ and \eqref{V1_phi0} follows from Lemma
\ref{lem_U}.

{\b Statements} (ii) and (iii) are corollaries from  \citet[Theorem
5.6]{feinberg2012reduction}.\hfill \Halmos \endproof

By  Theorem \ref{Th0}(iii), we consider only stationary policies
in the remaining part of the paper. Define
$V_{\alpha}^1(i,\delta)$ and $V_{\alpha}^0(i,\delta)$ {\b as}
\begin{align}
&V_{\alpha}^1(i,\delta)=(1-\delta)s_1+\displaystyle\frac{hi+c}{{\alpha+\lambda+i\mu }}
+\frac{\lambda}{\alpha+\lambda+i\mu }V_{\alpha}(i+1,1)+\frac{i\mu}{\alpha+\lambda+i\mu }V_{\alpha}(i-1,1),\notag\\
&V_{\alpha}^0(i,\delta)=\delta s_0+\displaystyle\frac{hi}{\alpha+\lambda}+\frac{\lambda}{\alpha+\lambda}V_{\alpha}(i+1,0).
\end{align}
%
%

\subsection{\uppercase{Full-Service Policies}}\label{sec_full}
 The class of policies that never switch the
running system off is the class of all policies for the case when
all the action sets $A(i,1)$ are  reduced to the set $\{1\}$ .
This is a sub-model of our original model defined by \eqref{optV}
with the action sets $A(i,1)$ reduced to $\{1\}$ for all $i =
0,1,\ldots$. Let $U_{\alpha}(i,\delta)$, $i = 0,1,\ldots$, be the
optimal {\b expected}  total discounted cost under the policies
that  never switch the system off.  

\begin{theorem}\label{Th1}  For any $\alpha>0$ the following statements
hold: 
\begin{hitemize}
\item[(i)]  for all $i=0,1,\ldots$,
$$
U_{\alpha}(i,1)=\displaystyle\frac{hi}{\mu+\alpha}+\frac{h\lambda}{\alpha(\mu+\alpha)}+\frac{c}{\alpha};
$$
\item[(ii)]   for all $i=0,1,\ldots$, the value function
$U_{\alpha}(i,0)$ satisfies the optimality equation
\begin{align}\label{DCOE_U}
U_{\alpha}(i,0)
= \min &\left\{
s_1+\displaystyle\frac{hi+c}{{\alpha+\lambda+i\mu }}
+\frac{\lambda}{\alpha+\lambda+i\mu }U_{\alpha}(i+1,1)
+\frac{i\mu}{\alpha+\lambda+i\mu }U_{\alpha}(i-1,1),\right.\notag\\
&\left.\
\frac{hi}{\alpha+\lambda}+\frac{\lambda}{\alpha+\lambda}U_{\alpha}(i+1,0)\right\}.
\end{align}
\end{hitemize}
\end{theorem}
\proof{Proof.} (i) For a policy $\phi$  that never switches the
running system off, $U_{\alpha}(i,1)=V_{\alpha}^{\phi}(i,1)$, and
the rest follows from Lemma \ref{lem_U}.

(ii)  Since $U_{\alpha}(i,0)$ is the optimal discounted cost for
the sub-model of the original MDP, it
 satisfies the discount-optimality equation. 
  Thus, \eqref{DCOE_U} follows from \eqref{DCOE}.
\hfill\Halmos\endproof

\begin{definition}\label{def_fullpol}
For an integer $n\ge0 $, a policy is called an $n$-full-service policy if it
never switches the running system off and switches the inactive system on
if and only if there are $n$ or more customers in the system. In
particular, the 0-full-service policy switches the system on at time
$0$, if it is off, and always keeps it on. A policy is called
a full-service policy if and only if it is an $n$-full-service policy for some $n\ge 0$.
\end{definition}

The following theorem implies that an $n$-full-service policy is discount-optimal within the class of policies that never switch the running system off.
\begin{theorem}\label{Th2}
A policy $\phi$ is discount-optimal within the class of the
policies that never switch off the running system if and only if
for all $i=0,1,\ldots$,
\[
\phi(i,0)=
\begin{cases}
1,  &\text{if $i>A(\alpha)$},\\
0,  &\text{if $i<A(\alpha)$},
\end{cases}
\]
where
\begin{equation}\label{A}
A(\alpha)=\displaystyle\frac{(\mu+\alpha)(c+\alpha s_1)}{h\mu}.
\end{equation}
\end{theorem}

Before proving Theorem \ref{Th2}, we introduce the definition of passive policies and some lemmas. In particular, the passive policy never changes the {\b status} of the system.

\begin{definition}\label{def_passive}
The policy $\varphi$,  with $\varphi(i,\delta)=\delta$ for all
$i=0,1,\ldots$ and {\b for }all $\delta=0,1$, is called passive.
\end{definition}

\begin{lemma}\label{lem_neveron}
For any $\alpha>0$, the passive policy $\varphi$ is not discount-optimal
within the class of policies that never switch off the running
system. Furthermore, $V_{\alpha}^{\varphi}(i,0)> U_{\alpha}(i,0)$
for all $i=0,1,\ldots$\ .
\end{lemma}
\proof{Proof.}
For the passive policy
$\varphi$,
 \begin{align*}
V_{\alpha}^{\varphi}(i,0)=\sum_{k=0}^{\infty}\left(\frac{\lambda}{\lambda+\alpha}\right)^k\frac{h(i+k)}{\lambda+\alpha}
=\frac{hi}{\alpha}+\frac{h\lambda}{\alpha^2}.
\end{align*}
For the policy $\phi$ that always runs the system,
\begin{equation}\label{eforfull} V_{\alpha}^{\phi}(i,0)=
s_1+\frac{hi}{\mu+\alpha}+\frac{h\lambda}{\alpha(\mu+\alpha)}+\frac{c}{\alpha}.
\end{equation} Thus
\begin{align*}
V_{\alpha}^{\varphi}(i,0)-V_{\alpha}^{\phi}(i,0)
&=\left(\frac{hi}{\alpha}+\frac{h\lambda}{\alpha^2}\right)-
\left(s_1+\frac{hi}{\mu+\alpha}+\frac{h\lambda}{\alpha(\mu+\alpha)}+\frac{c}{\alpha}\right)
=\frac{hi\mu}{\alpha(\mu+\alpha)}+\frac{h{\b \lambda\mu}}{\alpha^2({\b \mu}+\alpha)}-s_1-\frac{c}{\alpha}> 0,
\end{align*}
when $i$ is large enough. Let $i^*$ be the smallest natural integer such that the last inequality holds with $i=i^*$.
Let the initial state be $(i,0)$ with $i<i^*.$  Consider a policy
$\pi$ that keeps the system off in states $(j,0)$, $j<i^*,$ and
switches to a discount-optimal policy, when the number of
customers in the system reaches $i^*$.  Then
$V_\alpha^\varphi(i,0)>V_\alpha^\pi(i,0)\ge U_\alpha(i,0),$ where
the first inequality holds
because, before the process hits the state $(i^*,0)$, the policies $\varphi$ and $\pi$ coincide, and, after the process hits the state $(i^*,0)$, the policy $\pi$, which starting from that event coincides with $\phi$, incurs lower total discounting costs than the passive policy $\varphi$.
 \hfill\Halmos\endproof

\begin{lemma}\label{lem_Aalpha}
Let $\psi$ be the policy that switches the system on at time $0$ and keeps it on forever, and $\pi$ be the policy that waits for one arrival and then switches the system on and keeps it on forever. Then
\[\left\{\begin{array}{l}
V_{\alpha}^{\pi}(i,0)>V_{\alpha}^{\psi}(i,0),  \text{ if $i>A(\alpha)$},\\
V_{\alpha}^{\pi}(i,0)<V_{\alpha}^{\psi}(i,0),  \text{ if $i<A(\alpha)$},\\
V_{\alpha}^{\pi}(i,0)=V_{\alpha}^{\psi}(i,0),  \text{ if $i=A(\alpha)$},
\end{array}\right.
\]
where $A(\alpha)$ is as in \eqref{A}.
\end{lemma}
\proof{Proof.}
\begin{align*}
&V_{\alpha}^{\pi}(i,0)-V_{\alpha}^{\psi}(i,0)
=\left(\frac{hi}{\lambda+\alpha}+\frac{\lambda}{\lambda+\alpha}\left( s_1+U_{\alpha}(i+1,1) \right) \right)
 -\left( s_1+U_{\alpha}(i,1)\right)\notag\\
=&\left[\frac{hi}{\lambda+\alpha}
+\frac{\lambda}{\lambda+\alpha}\left(s_1+\frac{h(i+1)}{\mu+\alpha}+\frac{h\lambda}{\alpha(\mu+\alpha)}+\frac{c}{\alpha}\right)\right]
-\left[s_1+\frac{hi}{\mu+\alpha}+\frac{h\lambda}{\alpha(\mu+\alpha)}+\frac{c}{\alpha}\right]\notag\\
=&\frac{hi}{\lambda+\alpha}\frac{\mu}{\mu+\alpha}-\frac{\alpha}{\lambda+\alpha}\left(s_1+\frac{c}{\alpha}\right)
=\frac{h\mu}{(\lambda+\alpha)(\mu+\alpha)}\left(i-A(\alpha)\right),
\end{align*}
where the second equality holds in view of Theorem~\ref{Th1}(i)
and the rest is straightforward. \hfill\Halmos\endproof
%
\proof{Proof of
Theorem~\ref{Th2}.} Let $\phi$ be a stationary discount-optimal policy
within the class of the policies that never switch off the running
system. Let $\psi$ be the policy that switches the system on at time
$0$ and keeps it on forever, and $\pi$ be the policy that waits
for one arrival and then switches the system on and keeps it on
forever. By \eqref{DCOE_U},
\begin{align}\label{Th2_DCOE}
V_{\alpha}^{\phi}(i,0)=\min\left\{s_1+U_{\alpha}(i,1),\ \frac{hi}{\lambda+\alpha}+\frac{\lambda}{\lambda+\alpha}U_{\alpha}(i+1,0)\right\}.
\end{align}
We show that if $i>A(\alpha)$, then $\phi(i,0)=1$. {\b Assume} $\phi(i,0)=0$ for some $i>A(\alpha)$. By Lemma \ref{lem_neveron}, $\phi(j,0)=1$ for some $j>i$. Thus, there exists an $i^*\ge i$ such that $\phi(i^*,0)=0$ and $\phi(i^*+1,0)=1$. This implies that $V_{\alpha}^{\psi}(i^*,0)\ge V_{\alpha}^{\pi}(i^*,0)$, where $i^*>A(\alpha)$. By Lemma \ref{lem_Aalpha}, this is a contradiction. Thus $\phi(i,0)=1$ for all $i>A(\alpha)$.

If $i<A(\alpha)$, then Lemma \ref{lem_Aalpha} implies $V_{\alpha}^{\pi}(i,0)<V_{\alpha}^{\psi}(i,0)$. Thus $\phi(i,0)=0$ for all $i<A(\alpha)$.

Let $A(\alpha)$ be an integer and $i=A(\alpha)$.  In this case, Lemma \ref{lem_Aalpha} implies $V_{\alpha}^{\psi}(i,0)=V_{\alpha}^{\pi}(i,0)$.  From \eqref{DCOE_U},
$V_{\alpha}^{\psi}(i,0)=V_{\alpha}^{\pi}(i,0)=U_{\alpha}(i,0)=\min\left\{V_{\alpha}^{\psi}(i,0),V_{\alpha}^{\pi}(i,1)\right\}$. Thus $\phi(i,0)=1$ or $\phi(i,0)=0$. \hfill\Halmos\endproof

\begin{corollary}\label{cor_Th2}
Let
\begin{equation}\label{n}
n_{\alpha}=\lceil A(\alpha)\rceil,
\end{equation}
where $A(\alpha)$ is as in \eqref{A}. Then the following
statements hold: 

(i) if $A(\alpha)$ is not {\b an} integer then the $n_\alpha$-full
service policy is the unique stationary discount-optimal policy
within the class of policies that never switch the running system
off;

(ii) if $A(\alpha)$ is {\b an} integer then there are exactly two
stationary discount-optimal policies  within the class of policies
that never switch the running system off, and these policies are
$n_\alpha$- and $(n_\alpha+1)$-full-service policies;

(iii)
\begin{align}\label{Vbnd0}
U_{\alpha}(i,0)
=
\begin{cases}
\displaystyle\sum_{k=0}^{n_{\alpha}-i-1}\left(\displaystyle\frac{\lambda}{\lambda+\alpha}\right)^k\displaystyle\frac{h(i+k)}{\lambda+\alpha}
+\left(\frac{\lambda}{\lambda+\alpha}\right)^{n_{\alpha}-i}\left(s_1+\frac{hn_{\alpha}}{\mu+\alpha}
+\frac{h\lambda}{\alpha(\mu+\alpha)}+\frac{c}{\alpha}\right), &\text{ if $i<n_{\alpha}$},\\
s_1+\displaystyle\frac{hi}{\mu+\alpha}+\frac{h\lambda}{\alpha(\mu+\alpha)}+\frac{c}{\alpha}, &\text{ if $i\ge n_{\alpha}$}.
\end{cases}
\end{align}
\end{corollary}
\proof{Proof.}  Statements (i) and (ii) follow directly from
Theorem \ref{Th2} {\b and Definition~\ref{def_fullpol}}. Statements (i) and (ii) imply that
$V_\alpha^\phi=U_\alpha$, where $\phi$ is the $n_\alpha$-full
service policy. The first line of (\ref{Vbnd0}) is the
discounted cost to move from state $(i,0)$ to state
$(n_\alpha,0)$, when the system is off, plus the discounted cost
$U_\alpha(n_\alpha,0).$  The second line of (\ref{Vbnd0}) follows from
(\ref{eforfull}).
 \hfill\Halmos\endproof

\begin{corollary} \label{cor_Th2B} Let
 $n=\lfloor\displaystyle\frac{c}{h}+1\rfloor$.
 Then there exists
 $\alpha^*>0$ such that the  $n$-full-service policy is discount-optimal within
the class of the policies that never switch the running system off
for all discount rates $\alpha\in (0,\alpha^*].$

\end{corollary}
\proof{Proof.}
In view of \eqref{A}, the function $A(\alpha)$ is strictly monotone when $\alpha>0$. In addition, $A(\alpha)\searrow\displaystyle\frac{c}{h}$ when $\alpha\searrow 0$. This implies that $n_\alpha=n$ for all $\alpha\in(0,\alpha^{*}]$, where $\alpha^{*}$ can be found by solving the quadratic inequality
$
A(\alpha)\le n.
$
The rest follows from  Corollary~\ref{cor_Th2} (i) and (ii).
\hfill\Halmos\endproof

\subsection{\b\uppercase{Properties of Discount-Optimal Policies and Reduction to a Problem with a Finite State Space}}\label{sec_red}
{\b This subsection introduces the properties of the
discount-optimal policies formulated in Lemma \ref{lem_1} and
Lemma \ref{lem_Nbdd}, describes the inequalities between the major
thresholds in Lemma \ref{lem_Mbdd_old} that lead to the reduction
of the original infinite-state problem to a finite state problem.
This reduction essentially follows from  Corollary~\ref{Cor 3new}.
 Certain structural properties of discount-optimal policies are
described in Theorem~\ref{cor_disc}. }

Define
\begin{align}\label{M*}
M^*_{\alpha}=\begin{cases}\max\{i\ge 0:V_{\alpha}^0(i,1)\le
V_{\alpha}^1(i,1)\}, &\text{if $\{i\ge 0:V_{\alpha}^0(i,1)\le
V_{\alpha}^1(i,1)\}\ne\emptyset$},\\
-1, &\text{otherwise}.\end{cases}
\end{align}

\begin{lemma}\label{lem_1}
Let $\phi$ be a stationary discount-optimal policy. Then
$\phi(i,1)=1$ for $i\ge
\displaystyle\frac{h\lambda+c(\mu+\alpha)-s_0\alpha(\mu+\alpha)}{h\mu}$.
\end{lemma}
\proof{Proof.}
Let $\phi(i,1)=0$. Then $V_\alpha^\phi(i,1)> s_0+hi/\alpha$, since the number of customers in the system is always greater or equal than $i$ and after the first arrival it is greater than $i$. Observe that $V_\alpha^\phi(i,1)=V_\alpha(i,1)\le U_\alpha(i,1)$. From \eqref{U},
$$
s_0+\displaystyle\frac{hi}{\alpha}<\displaystyle\frac{hi}{\mu+\alpha}+\frac{h\lambda}{\alpha(\mu+\alpha)}+\frac{c}{\alpha}.
$$
This inequality implies
$i<\displaystyle\frac{h\lambda+c(\mu+\alpha)-s_0\alpha(\mu+\alpha)}{h\mu}$.
Thus, the opposite inequality implies
$\phi(i,1)=1.$\hfill\Halmos\endproof

\begin{corollary}\label{remark_M}
For all $\alpha >0$
\begin{equation}\label{eupestM}
M^*_\alpha \le
\frac{\lambda}{\mu}+\frac{(c+s_0\mu)^2}{4s_0h\mu}<\infty.
\end{equation}
\end{corollary}
\proof{Proof.}     According to Lemma~\ref{lem_1}, $M^*_\alpha\le
f(\alpha)$, where
$f(\alpha)=\displaystyle\frac{h\lambda+c(\mu+\alpha)-s_0\alpha(\mu+\alpha)}{h\mu}.$
For $\alpha>0$, the maximum of $f(\alpha)$ equals to the
expression on the right-hand side of
(\ref{eupestM}).\hfill\hfill\Halmos\endproof

\begin{lemma}\label{lem_Nbdd} Let $\phi$ be a stationary discount-optimal
policy. Then for any integer $j\ge 0$ there exists an integer
$i\ge j$ such that $\phi(i,0)=1.$
\end{lemma}
\proof{Proof.} If $\phi(i,0)={\b 0}$ for all $i\ge j$ then by
Lemma~\ref{lem_neveron}, $V_\alpha^\phi(j,0)>U_\alpha(j,0)\ge
V_\alpha(j,0).$  This contradicts the optimality of $\phi.$
\hfill\Halmos\endproof
%
%

Define
\begin{equation}\label{N*}
N^*_{\alpha}=\min\{i> M^*_\alpha:V_{\alpha}^1(i,0)\le
V_{\alpha}^0(i,0)\}.
\end{equation}
Lemma~\ref{lem_Nbdd} implies that $N_\alpha^*$ is well defined and
$N_\alpha^*<\infty$ for all $\alpha>0.$
{\b Before proving the relationship between $M^{*}_{\alpha}$ and
$N^{*}_{\alpha}$,  we introduce the following lemma.} 
\begin{lemma}\label{lem1}
The following properties hold for the function
$V_{\alpha}(i,\delta)$:
\begin{hitemize}
\item[(i)]  if $V_{\alpha}(i,0)=V_{\alpha}^1(i,0)$, then
$V_{\alpha}^1(i,1)<V_\alpha^0(i,1)$;
\item[(ii)] if $V_{\alpha}(i,1)=V_{\alpha}^0(i,1)$, then
$V_{\alpha}^0(i,0)<V_\alpha^1(i,0)$;
\item[(iii)] $ -s_1 \le
V_{\alpha}(i,1)-V_{\alpha}(i,0)\le s_0$.
\end{hitemize}
\end{lemma}
\proof{Proof.} (i) If $V_{\alpha}(i,0)=V_{\alpha}^1(i,0)$, then
$V_{\alpha}^1(i,0)\le V_{\alpha}^0(i,0)$. Hence
$V_{\alpha}(i,1)=V_{\alpha}(i,0)-s_1<
V_{\alpha}(i,0)+s_0=V_{\alpha}^0(i,1)$, where the inequality
follows from the assumption that $s_0+s_1>0.$  This implies
$V_{\alpha}^1(i,1)<V_\alpha^0(i,1)$.

(ii) If $V_{\alpha}(i,1)=V_{\alpha}^0(i,1)$, then
$V_{\alpha}^0(i,1)\le V_{\alpha}^1(i,1)$. Hence
$V_{\alpha}(i,0)=V_{\alpha}(i,1)-s_0<
V_{\alpha}(i,1)+s_1=V^1_{\alpha}(i,0)$.

(iii) $V_{\alpha}(i,0)\le s_1+V_{\alpha}(i,1)$  because
$V_{\alpha}(i,0)=\min\left\{s_1+V_{\alpha}(i,1),
V_{\alpha}^0(i,0)\right\}\le s_1+V_{\alpha}(i,1),$ and
$V_{\alpha}(i,1)\le s_0+V_{\alpha}(i,0)$  because
  $V_{\alpha}(i,1) =\min\left\{V_{\alpha}^1(i,1),
s_0+V_{\alpha}(i,0)\right\} \le s_0+V_{\alpha}(i,0).$
 \hfill\Halmos
\endproof

{\b The next Lemma shows the orders among $M_{\alpha}^{*}$,
$N_{\alpha}^{*}$ and $n_{\alpha}$.} {\b This leads to the
description of the properties of discount-optimal policies in
Corollary \ref{Cor 3new} that essentially reduces the problem to a
finite state space problem.}

\begin{lemma}\label{lem_Mbdd_old}
$M^*_{\alpha}<N^*_\alpha\le n_{\alpha}$ for all $\alpha>0$.
\end{lemma}
\proof{Proof.}  The definition~\eqref{N*} of $N^*_\alpha$ implies
that  $M^*_{\alpha}<N^*_\alpha$.  Thus, we need only to prove that
$N^*_\alpha\le n_{\alpha}$.

If $M_\alpha^*=-1$, according to \eqref{M*}, a discount-optimal
policy should never switch the running system system  off and
therefore $V_\alpha=U_\alpha$. In view of Corollary~\ref{cor_Th2},
$V_\alpha^0(i,0)<V_\alpha^1(i,0)$, when $i=0,\ldots,n_\alpha-1, $
and $V_\alpha^0(n_{\alpha},0)=V_\alpha^1(n_{\alpha},0).$   Thus,
in this case, $N^*_\alpha= n_{\alpha}$.

Let $M_\alpha^*\ge 0.$   Consider a stationary
discount-optimal policy $\varphi$ that switches the system on at state
  $(N_\alpha^*,0)$.  Such a policy exists in view of the
  definition of $N^*_\alpha$. It follows
   from the definition of $M^*_\alpha$ that
  $V^1_\alpha(i,1)<V^0_\alpha(i,1)$ for $i {\b >} M^{*}_\alpha$ .  Thus, the discount-optimal policy $\varphi$ always keeps running the active
  system
at states $(i, 1)$ when $i {\b >} M^*_\alpha$.
Observe that
  \begin{equation}\label{eq33} V^0_\alpha(N^*_\alpha-1,0)<V^1_\alpha(N^*_\alpha-1,0).\end{equation}
If $M^*_\alpha < N^*_\alpha-1$, \eqref{eq33} follows from the
definition of $N^*_\alpha.$ If $M^{*}_\alpha=N^*_\alpha-1$,
\eqref{eq33} follows from $V_\alpha^0(M_\alpha^*,1)\le
V_\alpha^1(M_\alpha^*,1)$ and from Lemma~\ref{lem1} (ii). Thus,
starting from the state $(N^*_\alpha-1,0)$, the discount-optimal
policy $\varphi$ waits until the next arrival, then switches the
system on and runs it until the number of customers in queue
becomes ${\b {M^*_\alpha\le }} N^*_\alpha-1$.  For $i=0,1,\ldots,$
let $F_{\alpha}^{1}(i)$ be the expected total discounted cost
incurred until the first time  $\theta(i)$ when the number of
customers in the system is $i$ and the system is running, if at
time 0 the system is off, there are $i$ customers in queue, and
the system is  switched on after the first arrival and is kept on
{\b as long as} the number of customers in system is greater than
$i$. Let $\theta=\theta(N^*_\alpha -1).$ Since $\varphi$ is the
discount-optimal policy,
$V_\alpha(N_\alpha^*-1,0)=F_{\alpha}^{1}(N_\alpha^*-1)+
[Ee^{-\alpha {\theta}}]V_\alpha(N_\alpha^*-1,1)$.

Let $\pi$ be a policy that switches the system on in state
$(N_\alpha^*-1,0)$ and then follows a discount-optimal policy.
Then, in view of \eqref{eq33}, the policy $\pi$ is not
discount-optimal at the initial state $(N^*_\alpha-1,0)$.  Thus,
$V^\pi_\alpha (N^*_\alpha-1,0)> V_\alpha(N_\alpha^*-1,0).$  Since
$V^\pi_\alpha (N^*_\alpha-1,0)=s_1+V_\alpha(N_\alpha^*-1,1)$,
\begin{equation*}
F_{\alpha}^{1}(N_\alpha^*-1)+ [Ee^{-\alpha
{\theta}}]V_\alpha(N_\alpha^*-1,1)<s_1+V_\alpha(N_\alpha^*-1,1),
\end{equation*}
and this is equivalent to
\begin{equation}\label{eq34}
\left(1-[Ee^{-\alpha
{\theta}}]\right)V_\alpha(N_\alpha^*-1,1)>F_{\alpha}^{1}(N_\alpha^*-1)-s_1.
\end{equation}

Assume that $n_\alpha<N_\alpha^*$.  Then $n_\alpha\le
N_\alpha^*-1$ and, in view of Theorem~\ref{Th2},
$\psi(N^*_\alpha-1,0)=1$ for a stationary discount-optimal policy
$\psi$ within the class of policies that never switches the system
off. Thus, $U_\alpha(N_\alpha^*-1,0)=V^{{\b
\psi}}_\alpha(N_\alpha^*-1,0)=s_1+U_\alpha(N_\alpha^*-1,1).$ In
addition, $U_\alpha(N_\alpha^*-1,0)\le V^{{\b
\varphi}}_\alpha(N_\alpha^*-1,0)=F_{\alpha}^{1}(N_\alpha^*-1)+[Ee^{-\alpha
{\theta}}]V_\alpha(N_\alpha^*-1,1).$ Thus,
\begin{equation}\label{eq35}
(1-[Ee^{-\alpha {\theta}}])U_\alpha(N_\alpha^*-1,1) \le
F_{\alpha}^{1}(i)-s_1.
\end{equation}
Observe that ${\theta}$ is greater than the time until the first
arrival, which has the positive expectation $\lambda^{-1}.$ Thus,
$[Ee^{-\alpha {\theta}}]<1$ and $U_\alpha(N_\alpha^*-1,1)\ge
V_\alpha(N_\alpha^*-1,1)$. \eqref{eq35} contradicts \eqref{eq34}.
Thus  $N_\alpha^*\le n_\alpha$. \hfill\Halmos\endproof

{
\begin{lemma}\label{lem_enewmew} For each $\alpha>0$,  the inequality $V_\alpha^1(i,0)\le V_\alpha^0(i,0)$  holds
when $i\ge n_\alpha.$
\end{lemma}
\proof{Proof.}  Fix any $\alpha>0$.  Consider two cases: case (i)
the best full-service policy is discount-optimal, and case (ii)
the best full-service policy is not discount-optimal.

Case (i).  According to Corollary~\ref{cor_Th2}, the
$n_\alpha$-full-service policy is discount-optimal.  This implies
that $V^1_\alpha(i,0)\le V^0_\alpha(i.0)$ for all $i\ge n_\alpha.$

Case (ii). {\b Let $\phi$ be a stationary discount-optimal
policy.} Assume that there exists an integer $j\ge n_\alpha$ such
that $\phi(j,0)=0$. Then, in view of Lemma~\ref{lem_Nbdd}, there
is $i\ge j$ such that $\phi(i,0)=0$ and  $\phi(i+1,0)=1.$ As shown
in Lemma~\ref{lem_Mbdd_old}, $n_\alpha>M^*_\alpha$ and therefore
$\phi(\ell,1)=1$ for all $\ell>M^*_\alpha$.  Thus,
$\phi(\ell,1)=1$ for all $\ell>i$.
 We have
\begin{align}\label{eq1_enewmew}
&V_\alpha^\phi(i,0) = F_{\alpha}^{1}(i) + [Ee^{-\alpha
\theta(i)}]V_\alpha(i,1)\le s_1+V_\alpha(i,1) \Rightarrow
F_{\alpha}^{1}(i)-s_1\le(1-[Ee^{-\alpha \theta(i)}])V_\alpha(i,1),
\end{align}
where the stopping time $\theta(i)$ and the expected total
discounted cost $F_{\alpha}^{1}(i)$ are defined in the proof of
Lemma~\ref{lem_Mbdd_old}.
  On the other
hand, since $i\ge n_\alpha$, under $n_\alpha$-full-service policy
$\pi$ we have
\begin{align}\label{eq2_enewmew}
V_\alpha^\pi(i,0)=s_1+U_\alpha(i,1)\le
F_{\alpha}^{1}(i)+[Ee^{-\alpha \theta(i)}]U_\alpha(i,1)
\Rightarrow (1-[Ee^{-\alpha \theta(i)}])U_\alpha(i,1)\le
F_{\alpha}^{1}(i)-s_1.
\end{align}
By \eqref{eq1_enewmew} and \eqref{eq2_enewmew}, we have
$U_\alpha(i,1)\le V_\alpha(i,1)$.  Since the best full-service
policy is not discount-optimal, $U_\alpha(i,1)> V_\alpha(i,1)$.  This
contradiction implies the correctness of the lemma.
\hfill\Halmos\endproof

\begin{corollary}\label{Cor 3new} Let $\alpha>0$ and $\alpha'\in (0,\alpha].$   For   a stationary
discount-optimal policy $\phi$ for the discount rate $\alpha'$,
consider the stationary policy $\phi'$,
\begin{equation}\label{eq:strarnew}
\phi'(i,\delta)=\begin{cases}
\phi(i,\delta), &\text{if $i<n_\alpha$},\\
 1, &\text{if $i\ge
n_\alpha$}.
\end{cases}
\end{equation}
Then the policy $\phi'$ is also discount-optimal for the discount
rate
$\alpha'$. 
\end{corollary}
\proof{Proof.} Let $\alpha'=\alpha$. By the definition \eqref{M*}
of $M^*_\alpha$, the inequality $V_\alpha^1(i,1)\le
V^0_\alpha(i,1)$ holds for all $i>M^*_\alpha$. By
Lemma~\ref{lem_enewmew} and by Corollary~\ref{cor_Th2},
$V_\alpha^1(i,0)\le V^0_\alpha(i,0)$ for all $i\ge n_\alpha$. In
view of Lemma~\ref{lem_Mbdd_old}, $M^*_\alpha<n_\alpha.$  Thus,
$V_\alpha^1(i,\delta){\b \le}V^0_\alpha(i,\delta)$ for all $i\ge
n_\alpha$ and for all $\delta=0,1$.  This implies the
discount-optimality of $\phi'$ for the discount rate $\alpha=\alpha'.$ Now let $\alpha'\in (0,\alpha).$
Since $\alpha>\alpha'>0$, then $n_{\alpha'}\le n_\alpha$, and 1 is an optimal decision for the discount rate $\alpha'$ at each state $(i,\delta)$ with $i\ge n_\alpha.$
 \hfill\Halmos\endproof

Corollary~\ref{Cor 3new} means that the system should be always
{\b run}, if there are $n_\alpha$ or more customers and the discount
rate is not greater than $\alpha.$  This essentially mean{\b s} that,
in order to find a discount-optimal policy for discount rates
$\alpha'\in (0,\alpha]$, the decision maker should find  such a policy only for a finite set of states $(i,\delta)$
with $i=0,1,\ldots,n_\alpha-1$ and $\delta=0,1.$
Thus, Lemma \ref{Cor 3new} reduces the original problem of optimization of the total discounted costs to a finite-state problem, and for every $\alpha >0$ this finite state set is the same  for all discount factors between 0 and $\alpha.$  The following theorem describes structural properties of a discount-optimal policy for a fixed discount factor.

\begin{theorem}\label{cor_disc} For each $\alpha> 0$, either the $n_\alpha$-full-service policy is discount-optimal, or there exists
a stationary discount-optimal policy $\phi_\alpha$ with the
following properties:
\begin{equation}\label{eqpropdiscopt}
\phi_\alpha(i,\delta)=\begin{cases} 1, &\text{if $i>M^*_\alpha$
and $\delta=1$},\\
1, &\text{if $i=N^*_\alpha$
and $\delta=0$},\\
1, &\text{{if $i\ge n_\alpha$ and $\delta=0$}},\\
0, &\text{if $i=M^*_\alpha$
and $\delta=1$},\\
0, &\text{if $M^*_\alpha\le i<N^*_\alpha$ and $\delta=0$}.
\end{cases}
\end{equation}
\end{theorem}
\proof{Proof.}  Consider a stationary discount-optimal policy
${\b \psi}$ for the discount rate $\alpha>0,$ and change it to
$\phi_\alpha$ according to \eqref{eqpropdiscopt} on the set of
states specified on the right-hand side of \eqref{eqpropdiscopt}.
The optimality of the new policy, denoted by $\phi_\alpha$,
follows from the definitions of $M^*_\alpha$ and $N^*_\alpha$, {and
from Corollary~\ref{Cor 3new}}.\Halmos\endproof

\section{\uppercase{The Existence and Structure of Average-Optimal Policies}}\label{s:avg}
{\b In this section we study the average cost criteria, prove the
existence of average-optimal policies and describe their
properties. }

\begin{definition}\label{def_mnpol}
For two nonnegative integers $M$ and $N$ with $N>M$, a stationary
policy is called an $(M,N)$-policy if
\begin{equation*}
\phi(i,\delta)=\begin{cases} 1, &\text{if $i>M$
and $\delta=1$},\\
1, &\text{if $i\ge N$
and $\delta=0$},\\
0, &\text{if $i\le M$
and $\delta=1$},\\
0, &\text{if $ i<N$ and $\delta=0$}.
\end{cases}
\end{equation*}
\end{definition}

\begin{theorem}\label{t:avopt}  There exists a stationary average-optimal
policy and, depending on the model parameters, either the $n$-full
service policy is average-optimal for $n=0,1,\ldots,$
 or an
$(M,N)$-policy is average-optimal for some $N>M\ge 0$ and $N\le n^*,$ where
\begin{equation}\label{n*}
n^*=\displaystyle{\lfloor\frac{c}{h}+1\rfloor}.
\end{equation}
In addition, the optimal average-cost value $v(i,\delta)$ is the
same for all initial states $(i,\delta)$; that is,
$v(i,\delta)=v.$
\end{theorem}
\proof{Proof.} We first prove that either the $n^*$-full-service
policy is average-optimal or an $(M,N)$-policy is average-optimal for some
$ N>M\ge 0$ and $N\le n^*$. For the initial {\b CTMDP}, consider a
sequence $\alpha_k\downarrow 0$ {\b as $k \to \infty$}. Let $\phi^k$
be a stationary discount-optimal policy for the discount rate
$\alpha_k$.   According to Theorem~\ref{cor_disc}, for each $k$
this policy can be selected either as {\b an} $n_{\alpha_k}$-full-service
policy or as a $\phi_{\alpha_k}$ policy satisfying
\eqref{eqpropdiscopt}. Since $n_{\alpha_k}\le
n_{\alpha_1}<(\mu+\alpha_1)(c+\alpha_1 s_1)/{h\mu}+1<\infty$ for
all $k=1,2,\ldots$, there exists a subsequence
{\b $\{{\alpha}_{k_\ell}\}$}, $\ell=1,2,\ldots,$ of the sequence
{\b $\{{\alpha}_k\}$}, $k=1,2,\ldots$ such that all the policies
$\phi^{k_\ell}=\phi,$ where $\phi$ is a stationary policy such
that either (i) the policy $\phi$ is an $n$-full-service policy
for some integer $n$ or (ii) the policy $\phi$ satisfies the
conditions on the right hand side of \eqref{eqpropdiscopt} with
the same $M^*_\alpha=M$ and $N^*_\alpha=N$ for
$\alpha=\alpha_{k_\ell}$.

Observe that the values of $v^\phi(i,\delta)$ do not depend on
the initial state $(i,\delta)$. Indeed, in case (i), when  the
{\b policy} $\phi$ is an $n^*$-full-service policy , the
stationary policy $\phi$ defines a Markov chain with a single
positive recurrent class $\left\{(i,1)\in Z:\, i=0,1,\ldots\right\}$, and all
the states in its complement  $\left\{(i,0)\in Z:\, i=0,1,\ldots\right\}$ are
transient. The same is true for case (ii) with the positive
recurrent class $Z^*=\{(i,1)\in Z:\, i=M,M+1,\ldots\}\cup
\left\{(i,0)\in Z:\, i=M,M+1,\ldots,N\right\} $ and with the set of transient
states $Z\setminus Z^*$.  In each case, the Markov chain leaves
the set of transient states in a finite expected amount of time
incurring a finite expected cost until the time the chain enters
the single positive recurrent class. Thus, in each case
$v^\phi(i,\delta)=v^\phi$ does not depend on $(i,\delta)$.

For all initial states $(i,\delta)$ and for an
arbitrary policy $\pi$, we have
\begin{align*}
v^\phi
=\lim_{t\to\infty}\ t^{-1}E^\phi_{(i,\delta)} C(t) {\b \le} \lim_{\alpha
\downarrow 0}\ \alpha V_\alpha^\phi(i,\delta)
\le \limsup_{\alpha \downarrow 0}\ \alpha V_\alpha^\pi(i,\delta)
\le \limsup_{t\to\infty}\ t^{-1}E^\pi_{(i,\delta)}
C(t)=v^\pi(i,\delta),
\end{align*}
where the first equality holds because of the definition of
average costs per unit time and the limit exists because {\b both}
$n^*$-full-service policy  and  $(M,N)$-policy define regenerative
processes, {\b the second and the last inequalities }follow from
the Tauberian theorem (see, e.g., \citet{korevaar2004tauberian}),
and the last {\b equality} is the definition of the average cost per
unit time. Since $\pi$ is an arbitrary policy, the policy $\phi$
is average-optimal. In addition, if $\alpha>0$ is sufficiently close to 0
then $n_\alpha=\lceil c/h\rceil$ if $c/h$ is not integer, and
$n_\alpha=c/h+1$, if $c/h$ is integer. This explains why
$n^*=\displaystyle{\lfloor\frac{c}{h}+1\rfloor}$ in
Theorem~\ref{t:avopt}. In conclusion, $v(i,\delta)=v$, since
$v^\phi(i,\delta)=v^\phi$.
In addition, if $n^*$-full-service policy $\phi$ is
average-optimal, and $\psi$ is an $n$-full-service policy for
$n=0,1,\ldots$, then $v^\psi=v^\phi$. 
\hfill\Halmos\endproof

\section{\uppercase{Computation of an Average-Optimal Policy}}\label{s:eg}

In this section, we show how an optimal policy can be computed via
Linear Programming.  According to Theorem~\ref{t:avopt}, there is an
optimal policy $\phi$ with $\phi(i,\delta)=1$ when $i\ge
n^*=\lfloor\displaystyle\frac{c}{h}+1\rfloor$. Thus, the goal is
to find the values of $\phi(i,\delta)$ when $i=0,1,\ldots,
{\b n^*-1}$ and $\delta=0,1.$  To do this, we truncate the state space
$Z$ to $Z^\prime=\{0,1,\ldots,n^*-1\}\times\{0,1\}.$  If the
action 1 is selected at state $(n^*-1,1)$, the system moves to the
state $(n^*-2,1)$, if the next change of the number of the
customers in the system is a departure and the system remains in
$(n^*-1,1)$, if an arrival takes place.  In the latter case, the
number of customers increases by one at the arrival time and then
it moves according to the random work until it hits the state
$(n^*-1,1)$ again.  Thus the system can jump from the state
$(n^*-1,1)$ to itself and therefore it cannot be described as a
CTMDP.  However, it can be described as a Semi-Markov Decision
Process (SMDP); see \citet[Chapter 5]{mine1970markovian} and \citet[Chapter 11]{puterman1994markov}.


We formulate our problem as an SMDP with the state set $Z^\prime$
and {\b the} action set $A(z)=A=\{0,1\}.$  If an action $a$ is selected at
state $z\in Z^\prime$, the system {\b spends} an average time
$\tau^\prime$ in this state until it moves to the next state
$z^\prime\in Z^\prime$ with the probability $p(z^\prime|z,a)$.
During this time the expected cost $C^{\prime}(z,a)$ is incurred.  {\b  For} ${\b z}=(i,\delta)$ with $i=0,1,\ldots, n^*-2$ and $\delta=0,1$,
these characteristics are the same as for the original CTMDP {\b and are given by}
\begin{align}\label{p}
p(z'|z,a)=
\begin{cases}
1, &\text{if $a=0,\ z^\prime=(i+1,0)$,}\\
\displaystyle\frac{\lambda}{\lambda+i\mu}, &\text{if $a=1,\
z^\prime=(i+1,1)$,}\\ \displaystyle\frac{i\mu}{\lambda+i\mu},
&\text{if $a=1,\ z^\prime=(i-1,1)$,}\\ 0, &\text{otherwise},
\end{cases}
\end{align}
\begin{align}\label{tau}
\tau' ((i,\delta),a)=
\begin{cases}
\displaystyle\frac{1}{\lambda}, &\text{if $a=0,$}\\
\displaystyle\frac{1}{\lambda+i\mu} &\text{if $a=1,$}
\end{cases}
\end{align}
and  $C' ((i,\delta),a)=|a-\delta|s_{a}+(hi+ac)\tau'
((i,\delta),a)$.  The transition probabilities in states
$(n^*-1,\delta)$ with $\delta=0,1$ are defined by
$p((n^*-2,1)|(n^*-1,\delta),1)=(n^*-1)\mu/(\lambda+(n^*-1)\mu)$,
$p((n^*-1,1)|(n^*-1,\delta),1)=\lambda/(\lambda+(n^*-1)\mu)$,
and
$p((n^*-1,1)|(n^*-1,\delta),0)=1.$ In the {\b last} case, the number
of customers increases by 1 to $n^*$, the system switches on, and
eventually the number of customers becomes $n^*-1$.

Let $T_i$ be the expected time between an arrival {\b seeing} $i$
customers in an $M/M/\infty$ queue and the next time when a
departure leaves $i$ customers behind, $i=0,1,\ldots\ .$ Applying
the memoryless property of the exponential distribution,
$T_{i}=B_{i+1}-B_{i}$, where $B_i$ is the expected busy period for
$M/M/\infty$ starting with $i$ customers in the system and $B_0=
0$. By formula (34b) in \citet{browne1995parallel},
\begin{align}\label{Tbusy'}
B_i
=\frac{1}{\lambda}\left(e^\rho-1+\sum_{k=1}^{i-1}\frac{k!}{\rho^k}\left(e^\rho-\sum_{j=0}^k\frac{\rho^j}{j!}\right)\right),
\end{align}
where $\rho=\displaystyle\frac{\lambda}{\mu}$. Thus $
T_{n^*-1}=B_{n^*}-B_{n^*-1}=\displaystyle\frac{1}{\lambda}\sum_{k=0}^\infty\frac{\rho^{k+1}}{n^*(n^*+1)\ldots(n^*+k)}.$

The expected time $\tau^\prime((n^*-1,\delta),1)$, where
$\delta=0,1$, is the expected time until the next arrival plus
$T_{n^*-1}$, if the next event is an arrival. Thus,
$\tau^\prime((n^*-1,\delta),1) =
\displaystyle\frac{\lambda}{\lambda+(n^{*}-1)\mu}\left(\frac{1}{\lambda}+T_{n^{*}-1}\right),$
$\delta=0,1.$ In addition $\tau^\prime((n^*-1,\delta),0) =
\displaystyle\frac{1}{\lambda}+T_{n^{*}-1},$
$\delta=0,1$.

To compute the one-step cost  $C^{\prime} ((n^*-1,1),1)$, we
define $m_{i}$ as the average number of visits to state $(i,1)$
starting from state $(n^{*}-1,1)$ and before revisiting state
$(n^{*}-1,1)$, $i=n^{*}-1,n^{*},\ldots\ $. And define $m_{i,i+1}$
as the expected number of jumps from $(i,1)$ to $(i+1,1)$,
$i=n^{*}-1,n^{*},\ldots$, and $m_{i,i-1}$  as the expected number
of jumps from $(i,1)$ to $(i-1,1)$, $i=n^{*},n^{*}+1,\ldots\ $.
Then $m_{i,i+1}=\displaystyle\frac{\lambda}{\lambda +i\mu}m_i$,
$m_{i,i-1}=\displaystyle\frac{i\mu}{\lambda +i\mu}m_i$ and
$m_{i,i+1}=m_{i+1,i}$. Since $m_{n^{*}-1}=1$,
\begin{align}\label{m}
m_{i}=\displaystyle\prod_{j=0}^{i-n^{*}}\displaystyle\frac{\lambda}{\lambda+(n^{*}-1+j)\mu}\frac{\lambda+(n^{*}+j)\mu}{(n^{*}+j)\mu},
\qquad i=n^{*},n^{*}+1,\ldots\ .
\end{align}
Thus,
\begin{align*}
C^{\prime}((n^*-1,1),1)&=\sum_{i=n^{*}-1}^{\infty}m_{i}C((i,1),1)
=\sum_{i=n^{*}-1}^{\infty}m_{i}\frac{hi+c}{\lambda+i\mu},
\end{align*}
where $C((i,1),1)=\displaystyle\frac{hi+c}{\lambda+i\mu}$,
$i=n^{*}-1,n^{*},\ldots$ is the cost incurred in state $(i,1)$
under action 1  for the original state space model; see Section
\ref{sec_OE}. The one-step cost $C^{\prime} ((n^*-1,0),1)=s_{1}+C'
((n^*-1,1),1)$.

Let $C_{n^*}$ be the total cost incurred in $M/M/\infty$ until the
number of customers becomes $(n^*-1)$ if at time 0 there are $n^*$
customers in the system and the system is running.  Then
\[C^{\prime}((n^*-1,1),1)=\displaystyle\frac{h(n^*-1)+c}{\lambda+(n^*-1)\mu}+\displaystyle\frac{\lambda}{\lambda+(n^*-1)\mu}C_{n^*},\]
and this implies
\[C_{n^*}=\left(1+\frac{(n^*-1)\mu}{\lambda}\right)C^{\prime}((n^*-1,1),1)-\frac{h(n^*-1)+c}{\lambda}.
\]
We also have
 $C^{\prime}
((n^*-1,0),0)=\displaystyle\frac{h(n^{*}-1)}{\lambda}+s_{1}+C_{n^*}
$, {\b $C^{\prime}
((n^*-1,0),1)=s_1 + C^{\prime}((n^*-1,1),1)
$}, and $C' ((n^*-1,1),0)=s_{0}+C' ((n^*-1,0),0)$.



With the definitions of the transition mechanisms, sojourn times,
and one-step costs for the SMDP, now we formulate the LP according
to Section 5.5 in \citet{mine1970markovian} or Theorem 11.4.2 and
formula 11.4.17 in \citet{puterman1994markov} as
\begin{equation}\label{lp_dual}
\begin{aligned}
\text{Minimize} \sum_{z\in Z^\prime}\sum_{a\in A}C' (z,a)x_{z,a}&\\
\text{s.t. \quad } \sum_{a\in A(z)}x_{z,a}-\sum_{z'\in Z^\prime}\sum_{a\in A(z)}{\b p(z|z^{\prime},a)x_{z,a}}&=0, \ z\in Z^\prime,\\
\sum_{z\in Z^\prime}\sum_{a\in A(z)}\tau' (z,a)x_{z,a}&=1,\\
x_{z,a}\ge 0, \ z\in Z^\prime ,\ a\in A.&
\end{aligned}
\end{equation}

Let $x^*$ be the optimal basic solution of \eqref{lp_dual}.
{\color{black} According to general results on SMDPs {\b in} \citet[Section
III]{denardo1970linear},} for each $z\in Z^\prime$, there exists at most
one ${\b a \in \{0,1\}}$ such that $x^*_{z,a}>0$. If $x^*_{z,a}>0$, then
for the average-optimal policy ${\b \phi,\ \phi(z)}=a$, for $a=0,1$. If
$x^*_{z,0}=x^*_{z,1}=0$, then ${\b \phi}(z)$ can be either 0 or 1.
For our problem, Theorem~\ref{Th_lp} explains how
$x^*{\b\vcentcolon=\{x^{*}_{z,a}:z\in Z', a\in A\}}$ can be used to construct a stationary average-optimal policy
$\phi$ with the properties stated in Theorem~\ref{t:avopt}.
{\color{black}
\begin{theorem}\label{Th_lp}
For an optimal basic solution $x^*$ of \eqref{lp_dual}, the
following statements hold: 
\begin{hitemize}
\item[(i)]if $x^*_{(0,1),1}>0$, then any $n$-full-service policy
is average-optimal, $n=0,1,\ldots;$ \item[(ii)]If
$x^*_{(0,1),0}>0$, then the $(0,N)$-policy is average-optimal with
\begin{align}\label{N_lp}
N=\begin{cases}
n^*,&\text{ if }\min\{i=1,\ldots,n^*-1:\,x^*_{(i,0),1}>0\}=\emptyset;\\
\min\{i=1,\ldots,n^*-1:\,x^*_{(i,0),1}>0\}, &\text{ if
}\min\{i=1,\ldots,n^*-1:\,x^*_{(i,0),1}>0\}\neq\emptyset;
\end{cases}
\end{align}
 \item[(iii)]if
$x^*_{(0,1),0}=x^*_{(0,1),1}=0$, then the $(M,N)$-policy is
average-optimal with
$M=\min\{i=1,\ldots,n^*-1:\,x^*_{(i,1),0}>0\}>0$ and $N$ being the
same as in \eqref{N_lp}.
\end{hitemize}
\end{theorem}
\proof{Proof.} Let $\phi^*$ be a stationary average-optimal policy
defined by the optimal basic solution $x^*$ of LP~\eqref{lp_dual}.
Since at most one of the values $\{x^*_{(0,1),0},x^*_{(0,1),1}\} $
is positive and they both are nonnegative, cases (i)--(iii) are
mutually exclusive and cover all the possibilities.

(i) If $x^*_{(0,1),1}>0$, then the state $(0,1)$ is recurrent
under the policy $\phi^*$ and $\phi^*(0,1)=1$.  Since the state
$(0,1)$ is recurrent and the system should be kept on in this
state, the policy $\phi^*$ always keeps the running system on.
The states corresponding to the inactive system are transient.
Thus, for any $n$-{\b full-service} policy $\phi$, $n=0,1,\ldots$, we have that
$v^\phi(j,0)=v^{\phi^*}(i,0)=v$ for all $i,j=1,2,\ldots\ .$ Thus,
any $n$-{\b full-service} policy is average-optimal.

(ii) If $x^*_{(0,1),0}>0$ then  the state $(0,1)$ is recurrent
under the policy $\phi^*$ and  $\phi^*(0,1)=0$. Since the state
$(0,1)$ is recurrent, the policy $\phi^*$ always keeps the running
system on as long as the system is nonempty.
By Lemma \ref{lem1} (ii), $\phi^*(0,0)=0$.
The first constraint in LP~\eqref{lp_dual}
implies that $x^*_{(1,0),0}+x^*_{(1,0),1}>0$. In general, if
$x^*_{(i,0),0}+x^*_{(i,0),1}>0$ for some $i=1,\ldots,n^*-1$, then
$\phi^*(j,0)=0$ if $x^{*}_{(j,0),1}=0$ for $j=0,\ldots,i-1$,  and $\phi^*(i,0)=1$ if $x^*_{(i,0),1}>0$.
Otherwise, if
$x^*_{(i,0),0}+x^*_{(i,0),1}=0$ for all $i=1,\ldots,n^*-1$,
$\phi^{*}(i,0)$ can be arbitrary and we define $\phi^*(i,0)=0$ for $i=0,1,\ldots,n^{*}-1$.  Thus,
formula \eqref{N_lp} defines the minimal number $N$ of customers
in the system, at which the inactive system should be switched on by
the average-optimal policy $\phi^*.$ We recall that the SMDP is defined
for the LP in the way that the system always starts on in state
$(n^*,0)$. Thus, the policy $\phi^*$ always keep running the
active system if the system is not empty, switches it off when the
system becomes empty, and switches on the inactive system when the
number of customers becomes $N$. If there are more than  $N$
customers when the system is inactive, the corresponding states
are transient. The defined $(0,N)$-policy starts the system
in all these states, and therefore it is average-optimal.

(iii) If $x^*_{(0,1),0}=x^*_{(0,1),1}=0$ then the state $(0,1)$ is
transient under the policy $\phi^*$.
In transient states  the average-optimal
policy $\phi^*$ can be defined arbitrary.  First observe that
$x^*_{(i,1),0}>0$ for some $i=1,\ldots,n^*-1$ and  therefore $M$
is well-defined in the theorem.  Indeed, if $x^*_{(i,1),0}=0$ for
all states $i=0,\ldots,n^*-1$, we can set $\phi^*(i,1)=1$ for all
these values of $i$.  This means that in the original Markov
chain, where the running system is always kept on when the number
of customers in the system is greater or equal than $ n^*$, the
system is always on. Since the birth-and-death for an $M/M/\infty$
system is positive recurrent, we have a contradiction.  Since the
state $(M,1)$ is recurrent for the Markov chain defined by the
policy $\phi^*$, this policy always keeps the running system on
when the number of the customers in the system is $M$ or more.
Since $x^*_{(i,\delta),a}=0 $ for $i<M$ and for all
$\delta,a=0,1$, we can define $\phi^*(i,\delta)$ arbitrarily when
$i<M.$  Let $\phi(i,\delta)=0$, when $i<M$ and $\delta=0,1.$
Similar to case (i), the policy $\phi^*$ prescribes to keep
inactive system off as long as the number of customers in the
{\b system is} less than {\b  $N$, switches it {\b on} when this number
becomes $N$, and it can be prescribed to switch the inactive system
on when the number of customers is greater than $N$, because all
such states are transient.  Thus, the defined $(M,N)$-policy is
optimal. \hfill\Halmos\endproof }

Similar to \eqref{lp_dual}, the LP can be formulated to find the
discount-optimal policy. However, this paper focuses on
average-optimality criteria, so we do not elaborate the LP for
discount-optimality here.

\section{\uppercase{Finding the Best $(0,N)$-Policy and Its Non-Optimality}}\label{s:no}

In this section we discuss how to compute the best $(0,N)$-policy and show that it may not be average-optimal. To do the latter, we consider an example.

Before providing the example, we show how to find the best $(0,N)$-policy.  This problem was studied by \citet{browne1995parallel} for the  $M/G/1$ queue without the running cost. Here we extend their solution to the case with running cost.
Let $\psi_{N}$ be a $(0,N)$-policy.
The average cost under $\psi_N$ can be found by formula (26) in \citet{browne1995parallel} by replacing the set up cost there with {\b the sum of switching costs and running costs} $s_{0}+s_{1}+c{\b B}_N$, where $B_{N}$ is as in \eqref{Tbusy'} or formula (34b) in \citet{browne1995parallel} .
This implies
\begin{equation}\label{v_0N}
v^{\psi_N}=hl_{N}+\frac{s_{0}+s_{1}+c{\b B}_{N}}{N/\lambda+B_{N}},
\end{equation}
 where $l_{N}$ is the expected number of customers in the system under $(0,N)$-policy. By formulae (22), (23) in \citet{browne1995parallel},
\begin{equation}\label{l_N}
l_{N}=\rho+\frac{N-1}{2}\frac{N}{N+\lambda B_{N}}.
\end{equation}
The optimal $N^{*}$ for the best $(0,N)$-policy is found by
\begin{equation}\label{best0N}
N^* = \argmin_N v^{\psi_N}.
\end{equation}
The following theorem extends Theorem 6 in \citet{browne1995parallel} to non-negative running cost $c\ge 0$.

\begin{theorem}\label{Th_Nbdd}
Let
\begin{equation}\label{Nbdd}
\tilde{N}=\min\left\{N\ge\frac{c}{h}\displaystyle\,:\,\frac{N(N+1)}{2\lambda}\ge\frac{s_0+s_1}{h}
 \right\},
\end{equation}
then for every $N\ge \tilde{N}$ we have that $v^{\psi_N}<v^{\psi_{N+1}}$, hence $\displaystyle \inf_{N\ge 1}v^{\psi_N}=\min_{1\le N\le\tilde{N}}v^{\psi_N}$.
\end{theorem}
\proof{Proof. }
To avoid notation conflict, let $b_N$ be $a_{N}$ defined as in formula (29) in \citet{browne1995parallel}. Note that $\displaystyle\frac{1}{\lambda}+T_N=b_{N-1}, N\ge 1$ and  $B_N+\displaystyle\frac{N}{\lambda}=\sum_{i=0}^{N-1}b_i$. By \eqref{v_0N} and \eqref{l_N}, we have for $(0,N)$-policy $\psi_N$ that
\begin{equation*}
v^{\psi_N}=h\left(\rho+\frac{N-1}{2}\frac{N/\lambda}{\sum_{i=0}^{N-1}b_i}\right)
 +\frac{s_0+s_1+c\left(\sum_{i=0}^{N-1}b_i-N/\lambda\right)}{\sum_{i=0}^{N-1}b_i}.
\end{equation*}
Thus $v^{\psi_N}<v^{\psi_{N+1}}$, if $hN-c> 0$ and
\begin{align*}
\left(\frac{h(N-1)}{2\lambda}+\frac{s_0+s_1}{N}-\frac{c}{\lambda}\right)/\left(\frac{hN}{\lambda}-\frac{c}{\lambda}\right)
<\displaystyle\frac{\sum_{i=0}^{N-1}b_i}{Nb_N}.
\end{align*}
Straightforward calculations show that, if $N\ge \tilde{N}$ in \eqref{Nbdd}, the left hand side $\le 1$ {\b in} the above inequality, and since the right hand side is always greater than 1 since $\left\{b_i: i\ge 0\right\}$ is decreasing, thus the result follows.
\hfill\Halmos\endproof

Theorem \ref{Th_Nbdd} implies that an average-optimal
$(0,N)$-policy can also be found by solving the LP~(\ref{lp_dual})
with the state space $Z''=\{(i,\delta):\ i=0,1,\ldots,
\tilde{N}-1,\ \delta=0,1\}$ and with the new action set
$A''(\cdot)$ defined as $A''(0,1)=\{0\}$, $A''(i,1)=\{1\}$ for
$i=1,\ldots,\tilde{N}-1$, and $A''(i,0)=\{0,1\}$ for
$i=1,\ldots,\tilde{N}-1$.

\begin{example}
We consider our model with arrival rate $\lambda=2$, service rate $\mu=1$ for each server, holding cost rate $h=1$, service cost rate $c=100$, and switching cost $s_0=s_1=100$. We implement it in \eqref{lp_dual} and run the LP with CPLEX in MatLab. We compute $n^*$ as $n^*=\lfloor\displaystyle\frac{c}{h}+1\rfloor=101$. Thus $Z^\prime=\{(i,\delta)\}$, with $i=0,1,\ldots,100$ and $\delta=0,1$. For the found solutions of \eqref{lp_dual},
$
x^*_{(i,0),1}>0, \text{ for }i=38;
x^*_{(i,0),0}>0, \text{ for }i=4,\ldots,38;
x^*_{(i,1),1}>0,\text{ for }i=5,\ldots,40;
x^*_{(i,1),0}>0,\text{ for }i=4;
$
and $x^*_{z,a}=0$ for all the other $z\in Z^{\prime},a\in A^{\prime}$. By Theorem \ref{Th_lp}, the average-optimal policy $\phi$ is $(M,N)$-policy with $M=4$ and $N=39$. The average cost of the $(4,39)$-policy is  $v^\phi\approx 43.39$.
The best $(0,N)$-policy can be found by solving \eqref{best0N}. Substituting \eqref{Tbusy'} and \eqref{l_N} in \eqref{v_0N}, we have
 $N^*=47$
and the corresponding average cost is  $v^{\psi_{N^*}}\approx 51.03 > v^\phi$.
%
\hfill\Halmos
\end{example}

\medskip
{\noindent\bf Acknowledgment\\}This research was partially
supported by NSF grants CMMI-0928490 and CMMI-1335296.  The
authors thank Manasa Mandava for useful comments.
%
%
%
%
%
%
%

\bibliographystyle{apalike2}
\bibliography{XBib}

\end{document}